\documentclass[12pt]{article}
\usepackage[russian,english]{babel}
\usepackage{amsmath,amsfonts,amssymb,amsthm,graphics,epsfig}

\newtheorem{thm}{Theorem}[section]
\newtheorem{lem}[thm]{Lemma}

\newtheorem{con}{Conjecture}
\theoremstyle{definition}
\newtheorem{defin}[thm]{Definition}

\newcommand { \ib }[1] {\textit{\textbf{#1}}}

\begin{document}
\renewcommand{\ib}{\mathbf}
\renewcommand{\proofname}{Proof}
\renewcommand{\phi}{\varphi}
\makeatletter \headsep 10 mm \footskip 10 mm
\renewcommand{\@evenhead}%

\title{The second Voronoi conjecture on parallelohedra for zonotopes\footnote{This work is financially supported by RFBR (projects 11-01-00633-a and 11-01-00735-a) and by grant of the President of the Russian Federation  {\Russian อุ}-5413.2010.1.}.}
\author{A. Garber\footnote{Lomonosov Moscow State University, Faculty of Mechanics and Mathematics, Department of Higher Geometry and Topology. Russia, 119991, Moscow, Vorob'evy gory, 1, A-1620.}\\
{\tt e-mail: alexeygarber@gmail.com}}

\maketitle

\begin{abstract}

We prove the second Voronoi conjecture on parallelohedra for zonotope. We show that for a given face-to-face tiling of $d$-dimensional Euclidean space into parallel copies of zonotope $Z$ there are $d$ vectors, connecting centers of zonotopes with common facet, that are basis of the correspondent lattice of the tiling.

\end{abstract}

\section{Parallelohedra}

\begin{defin}
A polytope $P\subseteq \mathbb{R}^d$ is called a {\it parallelohedron} if Euclidean space $\mathbb{R}^d$ can be tiled into non-overlapping parallel copies of $P.$
\end{defin}

There are several classical results in the theory of parallelohedra. Here we will mention two of them. The first one is the Minkowski theorem \cite{Min} that claims that parallelohedron $P$ is centrally symmetric, any facet of $P$ is centrally symmetric and projection of $P$ along any of its ridge (i.e. face of codimension $2$) is either parallelogram or centrally symmetric hexagon. The second one is the Venkov theorem that claims that these three conditions of Minkowski are sufficient for $P$ to be a parallelohedron.

In general case two parallelohedra with common boundary point can share only a part of faces of both, as in usual brickwork two bricks from consecutive horizontal layers share only half of facet.

\begin{defin}
If intersection of any two polytopes of the tiling $\mathcal{T}$ of $\mathbb{R}^d$ into non-overlapping polytopes is a face of both of them (this face can be empty) then the tiling $\mathcal{T}$ is called {\it face-to-face} tiling.
\end{defin}

McMullen showed \cite{McM} that if parallelohedron $P$ admits an arbitrary tiling of $\mathbb{R}^d$ into parallel copies then it admits also and a face-to-face tilings, namely McMullen proved that Minkowski conditions are necessary for parallelohedron $P$ with non face-to-face tiling. Further for a given parallelohedron $P$ we will consider only correspondent face-to-face tiling $\mathcal{T}(P).$ The uniqueness (up to translation) of such face-to-face tiling into copies of $P$ is evident. Centers of all polytopes of the tiling $\mathcal{T}(P)$ forms a $d$-dimensional lattice $\Lambda(P).$

Also for an arbitrary $d$-dimensional lattice $\Lambda$ we can construct a $d$-dimensional parallelohedron. For a fixed point $O\in\Lambda$ consider a polytope $P_\Lambda$ that consists of all points of $\mathbb{R}^d$ that are closer to $O$ than to any other point of $\Lambda.$ The constructed polytope $DV_\Lambda$ is called the {\it Dirichlet-Voronoi polytope} for $\Lambda$ and this polytope is a parallelohedron because for different lattice points corespondent polytopes differs only by translation and all such polytopes gives us a face-to-face tiling of $\mathbb{R}^d.$

The first conjecture of Voronoi states that all parallelohedra can be obtained from Dirichlet-Voronoi polytopes with affine transformations.

\begin{con}[G.Voronoi \cite{Vor}]\label{con1}
For any parallelohedron $P$ there exists a lattice $\Lambda$ and an affine transformation $\mathcal{A}$ such that the polytope $\mathcal{A}(P)$ is the Dirichlet-Voronoi polytope for $\Lambda.$
\end{con}

Consider the set $\mathcal{N}(P)$ of vectors connecting the center of a given parallelohedron $P$ with centers of all other parallelohedra of tiling $\mathcal{T}(P)$ that shares facets with $P.$ It is clear that $\mathcal{N}(P)$ generates the lattice $\Lambda(P)$ because any vector from $\Lambda(P)$ can be represented as a sum of several vectors from $\mathcal{N}(P).$

\begin{con}[G.Voronoi \cite{Vor}]\label{con2}
We can choose $d$ vectors from $\mathcal{N}(P)$ that forms a basis of $\Lambda(P).$
\end{con}

In this paper we will prove the second conjecture of Voronoi in the case of space-filling zonotopes.

\section{Delone tilings and Dirichlet-Voronoi tilings}

We can generalize the construction of Dirichlet-Voronoi polytope on the case of general positive definite quadratic form.

\begin{defin}
Let $\phi: \mathbb{R}^d \longrightarrow \mathbb{R}$ be a positive definite quadratic form in $d$-dimensional Euclidean space and let $O$ be a point of some fixed lattice $\Lambda.$ The {\it Dirichlet-Voronoi polytope $P_\phi(\Lambda)$} for $\Lambda$ with respect to form $\phi$ is a polytope that consist of all points $X$ of $\mathbb{R}^d$ such that value of $\phi$ on vector $\overrightarrow{OX}$ is not greater than value $\phi(\overrightarrow{O'X})$ for any other point $O'\in\Lambda.$
\end{defin}

The constructed polytope $P_\phi(\Lambda)$ is a parallelohedron too and the correspondent tiling is called a {\it Dirichlet-Voronoi tiling} $\mathcal{V}_\phi(\Lambda)$. If we will take $\phi$ equals to the usual metric form $\phi(\ib x)=\ib x^T\ib x$ then we will get the usual Dirichlet-Voronoi polytope $DV_\Lambda.$

Applying an affine transformation $\mathcal{A}$ to a given lattice $\Lambda$ and a polytope $P_\phi(\Lambda)$ with form $\phi(\ib x)=\ib x^T Q\ib x$ we will obtain lattice $\mathcal{A}\Lambda$ and  polytope $P_{\phi_\mathcal{A}}(\mathcal{A}\Lambda)$ with respect to quadratic form $\phi_\mathcal{A}(\ib x)=(\mathcal A^{-1}\ib x)^T Q (\mathcal A^{-1}\ib x).$ So to prove the first conjecture of Voronoi for a given polytope $P$ with lattice $\Lambda$ it is enough to show that there exist a quadratic form $\phi$ such that $P=P_\phi(\Lambda).$

For any Dirichlet-Voronoi tiling of $\mathbb{R}^d$ into polytopes $P_\phi(\Lambda)$ we can construct a dual Delone tiling. The Delone tiling $\mathcal{D}_\phi(\Lambda)$ is defined by lattice $\Lambda$ and positive definite quadratic form $\phi.$ For a given form $\phi$ consider an ellipsoid $\mathcal E_\phi$ with equation $\phi(\ib x)=1.$ Consider an arbitrary homothetic copy $\mathcal E_0$ of $\mathcal E_\phi$ such that $\mathcal E_0$ has some $d$-dimensional set of points from $\Lambda$ on the boundary $\partial \mathcal E$ but does not have points from $\Lambda$ inside $\mathcal E_0.$ This ellipsoid $\mathcal E_0$ defines a convex polytope inscribed in $\mathcal E$ with vertices from the lattice $\Lambda.$ The set of all such polytopes inscribed in ``empty'' ellipsoids we will call the {\it Delone tiling} $\mathcal{D}_\phi(\Lambda).$

For any vertex $O$ of the Dirichlet-Voronoi tiling values of quadratic form $\phi$ on vectors that connects $O$ with centers of polytopes from $\mathcal{V}_\phi(\Lambda)$ that meets at $O$ are equal so $O$ will be center of empty ellipsoid that defines one polytope from the Delone tiling $\mathcal{D}_\phi(\Lambda).$

\section{Zonotopes and dicings}

\begin{defin}
A polytope $P\subseteq \mathbb{R}^d$ is called a {\it zonotope} if $P$ can be represented as a Minkowski sum of finite number of segments. The Minkowski sum of segments defined by vectors of the set $V=\{\ib v_1,\ldots,\ib v_n\}$ we will denote by $Z(V).$ Vectors $\ib v_i$ is called the {\it zone vectors} of the zonotope $Z(V)$.

Equivalently zonotope can be defined as a projection of cube $C^n$ of some dimension $n \geq d.$
\end{defin}

Erdahl in his work \cite{Erd} proved the first conjecture of Voronoi for zonotopes. Later Deza and Grishukhin proved the first conjecture of Voronoi for zonotopes in terms of oriented matroids \cite{DG}. In this paper we will formulate and use several notions and theorems concerning zonotopes and correspondent Delone tilings and Dirichlet-Voronoi tilings from Erdahl's work \cite{Erd}.

\begin{defin}
Consider $n$ families of hyperplanes in $d$-dimensional Euclidean space. Assume that every family consists of parallel hyperplanes and slices $\mathbb R^d$ into layers of a constant width (this width can vary for different families). This set of $n$ families and the correspondent tiling of $\mathbb{R}^d$ is called {\it dicing} if there are $d$ linearly independent normals to these hyperplanes and every point that belongs to hyperplanes of $d$ families with independent normals also belongs to hyperplanes of all other families.
\end{defin}

Assume that the point with radius vector $\ib a_0$ belongs to hyperplanes of all families then for a given dicing we can choose normals $\ib d_1,\ldots,\ib d_n$ to families in such a way that all hyperplanes of dicing will be defined by equations $\ib d_i\cdot(\ib x-\ib a_0)=a$ for various integer $a.$ Also we can substitute any vector $\ib d_i$ by its negative $-\ib d_i.$

\begin{defin}
In that case we will denote our dicing by $\mathfrak{D}(\ib d_1,\ldots,\ib d_n)$ and the set $\mathcal{D}=\{\pm\ib d_1,\ldots,\pm\ib d_n\}$ we will call the {\it set of normals} of the dicing $\mathfrak{D}.$
\end{defin}

Consider an arbitrary vertex of a dicing $\mathfrak{D}$ (i.e. point of intersection of hyperplanes from $d$ families with independent normals) and consider all edges of $\mathfrak{D}$ that incident to this vertex. This set of edges will be centrally symmetric.

\begin{defin}
The mentioned set of edges $\mathcal{E}=\{\pm\ib e_1,\ldots,\pm\ib e_k\}$ is called the {\it edge set} of a dicing $\mathfrak{D}.$ It is clear that construction of the set $\mathcal{E}$ does not depend on the vertex of a dicing.

Also this construction immediately follows that all vertices of a dicing $\mathfrak{D}$ forms a lattice $\Lambda(\mathfrak{D})$ that we will call the {\it lattice of a dicing}.
\end{defin}

In particulary Erdahl prooved the following theorem on connection between sets $\mathcal{D}$ and $\mathcal{E}$ \cite[Theorem 3.1]{Erd}.

\begin{thm}[R. Erdahl]\label{edges}
Given set $\mathcal D$ of $n$ vectors with $d$ linearly independent vectors can be a set of normals for a dicing if and only if there exists a set $\mathcal E$ such that:

\begin{enumerate}
\item[{\rm (E1)}] Any pair of opposite vectors $\pm \ib e_i\in \mathcal{E}$ lies in some one-dimensional intersection $\ib d_{i_1}^\perp\cap\ldots\cap \ib d_{i_{d-1}}^\perp$ with independent vectors $d_{i_j}\in \mathcal D$ and conversely for any $d-1$ linearly independent vectors from $\mathcal D$ there is a correspondent pair of opposite vectors from $\mathcal E;$

\item[{\rm (E2)}] For any pair of vectors $\ib d\in \mathcal D$ and $\ib e\in \mathcal E$ the scalar product $\ib d^T\ib e$ equals to $0$ or $\pm 1.$
\end{enumerate}
\end{thm}

\begin{defin}
Matrices $\mathbf{D}$ and $\mathbf{E}$ with vectors from $\mathcal D$ and $\mathcal E$ written in columns are called {\it matrix of normal vectors} and {\it matrix of edge vectors} of the dicing $\mathfrak{D}$ respectively.
\end{defin}

Under the affine transformation of $\mathbb R^d$ with matrix $L$ matrices $\mathbf{D}$ and $\mathbf{E}$ changes into matrices $\mathbf{D}'=(L^{-1})^T\mathbf{D}$ and $\mathbf{E}'=L\mathbf{E}.$ Moreover, there exists an affine transformations such that entries of matrices $\mathbf{D}$ and $\mathbf{E}$ after this transformation will be only $0$ and $\pm 1$ \cite[Theorem 3.3]{Erd}.

\begin{thm}[R. Erdahl]
There exists an affine transformation that will give us a totally unimodular matrix $\mathcal D'$, i.e. any of its minor will be equal to $0$ or $\pm 1.$ Moreover this transformation can be chosen in such a way that both sets $\mathbf{D}$ and $\mathbf{E}$ will contain $d$ vectors of standard basis of $\mathbb R^d$, i.e. vectors $(1,0,\ldots, 0)^T,(0,1,0,\ldots, 0)^T,\ldots,(0,\ldots,0,1)^T.$ Also after this transformation all entries of matrices $\mathbf{D}$ and $\mathbf{E}$ will become $0$ or $\pm 1.$
\end{thm}

We will take only a half of vectors from the set $\mathcal{D}$ with no opposite vectors included. This new set $\mathcal{D}^+$ determines the same unique dicing $\mathfrak D.$ Consider a quadratic form
$$\phi(\ib x)=\sum_{\ib d\in \mathcal{D}^+}\omega_{\ib d}(\ib d\cdot \ib x)^2 $$
for some positive constants $\omega_{\ib d}.$ The correspondent Dirichlet-Voronoi polytope described in \cite[Theorem 4.3]{Erd}.

\begin{thm}[R. Erdahl]\label{zones}
The Dirichlet-Voronoi polytope for lattice $\Lambda(\mathfrak{D})$ with respect to quadratic form $\phi (\ib x)$ is a zonotope with zone vectors $$\ib z_{\ib d}=\left(\sum_{\ib d\in\mathcal{D}^+}\omega_{\ib d}\ib d\ib d^T\right)^{-1}\omega_{\ib d}\ib d, \quad \ib d\in \mathcal D^+.$$
\end{thm}

In the same work Erdahl proved that if $\phi(\ib x)$ is the standard Euclidean metrics then the zone vectors of the correspondent zonotope are $\ib z_{\ib d}=\omega_{\ib d}\ib d$ (\cite[Sect. 6, p. 442]{Erd}).

Also \cite[Theorem 1.2]{Erd} claims that Dirichlet-Voronoi polytope for a lattice is a zonotope if and only if the correspondent Delone tiling is a dicing.

\section{The second conjecture of Voronoi for zonotopes}

\begin{lem}\label{face}
Let $Z=Z(V)$ be a $d$-dimensional zonotope. Any facet of $Z$ is generated by some $(d-1)$-dimensional subset $U$ of $V$ and conversely any $(d-1)$-dimensional subset of $V$ generates facet $Z(U)$ of the zonotope $Z(V).$
\end{lem}

\begin{proof}
Let $\pi$ be a hyperplane of some facet of $Z.$ The zonotope $Z(V)=Z(\ib v_1,\ldots,\ib v_n)$ is a projection of a cube $C^n\subset\mathbb{R}^n$ onto space $\mathbb{R}^d$ along $(n-d)$-dimensional subspace $\psi.$ Consider the hyperplane $\pi \times \psi$ in the space $\mathbb R^n.$ This hyperplane is a supporting plane of the cube $C^n$ and hence it defines its face $F$. The face $F$ is a cube of some dimension and this face is generated by edges of $C^n$ that projects onto vectors of the set $V$ that are parallel to $\pi.$ Hence the face $F$ is projected into parallel copy of a zonotope $Z(U)$ for some $(d-1)$-dimensional subset $U$ of $V.$ The converse statement can be proven in the analogous way.
\end{proof}

\begin{thm}
The conjecture $\ref{con2}$ is true for space filling zonotopes.
\end{thm}

\begin{proof}
Assume that $Z$ is a $d$-dimensional space filling zonotope, i.e. zonotope that is also a parallelohedron. Then the first conjecture of Voronoi is true for $Z$ \cite[Theorem 1.1]{Erd}, i.e. there exists an affine transformation $\mathcal{A}$ such that the zonotope $\mathcal{A}Z$ is Dirichlet-Voronoi polytope of some lattice $\Lambda$ with respect to the usual Euclidean metrics as a quadratic form. The Delone tiling for $\Lambda$ with respect to Euclidean metrics is a dicing $\mathfrak{D}=\mathcal{D}(\pm\ib d_1,\ldots,\pm\ib d_n)$ so due to theorem \ref{zones} the zone vectors of the zonotope $\mathcal{A}Z$ can be written as $\omega_i \ib d_i.$

Any vector that connects centers $O_1$ and $O_2$ of copies of $\mathcal{A}Z$ with a joint facet $F$ is perpendicular to this facet because any point of $F$ is equidistant from $O_1$ and $O_2.$ By lemma \ref{face} there are $d-1$ vectors from the set $\left\{\omega_i\ib d_i\right\}_{i=1}^n$ that are parallel to $F$ so by theorem \ref{edges} there is a vector from $\mathcal{E}$ perpendicular to $F$. That means the the vector $\overrightarrow{O_1O_2}$ lies in $\mathcal{E}.$ The converse statement is also true, if we take any vector $\ib x$ from $\mathcal{E}$ then $d-1$ linearly independent vectors from $\mathcal{D}$ that are perpendicular to $\ib x$ will determine a facet of $Z.$ Therefore sets $\mathcal{E}$ and $\mathcal{N}(Z)$ coincides.

Consider a unimodular representation of the dicing $\mathfrak{D}.$ In this representation the set $\mathcal{E}$ contains the standard basis of the space $\mathbb{R}^d$ and all other vectors of $\mathcal{E}$ has integer coordinates. Hence in unimodular representation the set $\mathcal{N}(Z)$ contains a standard basis of $\mathbb{R}^d$ and all other vectors of this set has integer coordinates in this basis. So in unimodular representation the lattice $\Lambda(Z)$ coincides with $\mathbb{Z}^d$ and $\mathcal{N}(Z)$ contains its basis. The last property is invariant under linear transformations and changing of the basis in $\mathbb{R}^d$ and this proves the theorem.
\end{proof}


\begin{thebibliography}{20}
\bibitem{DG} M. Deza, V. Grishukhin, \emph{Once more about Voronoi's conjecture and space tiling zonotopes}. e-print, http://arxiv.org/abs/math/0203124, 2002.

\bibitem{Erd} R.~Erdahl, \emph{Zonotopes, Dicings, and Voronoi's Conjecture on Parallelohedra}, Eur. J. of Comb., Vol. 20, N. 6, 1999, pp. 527-549.

\bibitem{McM} P.McMullen, \emph{Convex bodies which tile space by translation}. Mathematika, vol. 27, n. 1, 1980, pp. 113-121.

\bibitem{Min} H.~Minkowski, \emph{Allgemeine Lehrs\"atze \"uber die convexen Polyeder}. G\"ott. Nachr., 1897, pp. 198-219.

\bibitem{Vor} G.~Voronoi, \emph{Nouvelles applications des param\`etres continus \`a la th\'eorie des formes quadratiques. Deuxi\`eme m\'emoire. Recherches sur les parall\'elo\`edres primitifs}. J. f\"ur Math., vol. 136, 1909, pp. 67-178.

\bibitem{Ven} B.A.~Venkov, \emph{About one class of Euclidean polytopes} (in Russian). Vestnik Leningr. Univ., ser. Math., Phys., Chem., 1954, vol. 9, pp. 11-31.

\end{thebibliography}
\end{document}

\renewcommand{\ib}{\mathbf}
\renewcommand{\proofname}{Доказательство}
\renewcommand{\phi}{\varphi}
\makeatletter \headsep 10 mm \footskip 10 mm
\renewcommand{\@evenhead}%

\title{Вторая гипотеза Вороного о параллелоэдрах для зонотопов\footnote{Работа выполнена при финансовой поддержке РФФИ (проекты 11-01-00633-a и 11-01-00735-a), а также при поддержке гранта Президента РФ НШ-5413.2010.1.}.}
\author{А.И. Гарбер\footnote{Московский государственный университет им. М.В.Ломоносова, механико-математический факультет, каф. высшей геометрии и топологии. Россия, 119991, Москва, Воробьевы горы, 1, А-1620.}\\
{\tt e-mail: alexeygarber@gmail.com}}

\maketitle

\begin{abstract}
В данной работе доказана вторая гипотеза Вороного о параллелоэдрах для зонотопов, то есть доказано, что для данного разбиения ``грань-в-грань'' $d$-мерного пространства на копии зонотопа $Z$ можно выбрать $d$ векторов, соединяющие центры зонотопов соседних по гиперграням, которые являются базисом соответствующей решетки разбиения.
\end{abstract}

\section{Параллелоэдры.}

\begin{defin}
Многогранник $P\subseteq \mathbb{R}^d$ называется {\it параллелоэдром} если пространство $\mathbb{R}^d$ можно разбить на параллельные копии $P$, не пересекающиеся по внутренним точкам.
\end{defin}

Классическими результатами в теории параллелоэдров являются теорема Минковского \cite{Min}, которая утверждает, что если многогранник $P$ является параллелоэдром, то он сам имеет центр симметрии, любая его гипергрань имеет центр симметрии и ортогональная проекция $P$ вдоль любой грани коразмерности $2$ есть параллелограмм или центрально-симметричный шестиугольник, и теорема Венкова \cite{Ven}, которая утверждает, что данные три условия являются достаточными для того, чтобы многогранник $P$ был параллелоэдром.

В общем случае два соседних параллелоэдра могут пересекаться по части грани каждого из них, например, в обычной ``кирпичной кладке'' два соседних кирпича из разных горизонтальных слоев пересекаются не по полной грани.

\begin{defin}
Если в разбиении $\mathcal{T}$ любые два параллелоэдра пересекаются по грани каждого из них (эта грань может быть и пустой), то соответствующее $\mathcal{T}$ называется разбиением {\it грань-в-грань}.
\end{defin}

МакМаллен \cite{McM} доказал, что если данный параллелоэдр $P$ допускает разбиение пространства $\mathbb{R}^d$ на параллельные копии произвольным образом, то он допускает и разбиение грань-в-грань, а именно МакМаллен доказал, что даже для параллелоэдра $P$ и разбиения не грань-в-грань три условия Минковского являются необходимыми. В дальнейшем для каждого параллелоэдра $P$ мы будем рассматривать только соответствующее ему разбиение $\mathcal{T}(P)$ грань-в-грань. Единственность такого разбиения очевидна. При этом центры многогранников разбиения $\mathcal{T}(P)$ образуют $d$-мерную решетку $\Lambda(P).$

Также по любой $d$-мерной решетке $\Lambda$ можно построить $d$-мерный параллелоэдр. Для фиксированной точки $O\in \Lambda$ рассмотрим многогранник $P_\Lambda,$ состоящий из всех точек пространства $\mathbb{R}^d,$ которые ближе к $O$ чем к остальным точкам $\Lambda.$ Данный многогранник называется {\it многогранником Дирихле-Вороного} и является параллелоэдром, так как многогранники для различных точек решетки отличаются параллельным переносом и дают разбиение ``грань-в-грань'' пространства $\mathbb{R}^d.$

Первая гипотеза Вороного утверждает, что все параллелоэдры могут быть получены с помощью аффинных преобразований из многогранников Дирихле-Вороного.

\begin{con}[Г.Вороной \cite{Vor}]\label{con1}
Для любого параллелоэдра $P$ существуют такие решетка $\Lambda$ и аффинное преобразование $\mathcal{A},$ что многогранник $\mathcal{A}(P)$ является многогранником Дирихле-Вороного для решетки $\Lambda.$
\end{con}

Пусть $\mathcal{N}(P)$ --- набор векторов соединяющих центр данного параллелоэдра $P$ с центрами его соседями по гиперграням в разбиении $\mathcal{T}(P).$ Очевидно, что множество $\mathcal{N}(P)$ порождает решетку $\Lambda(P)$ так как любой вектор, соединяющий центры двух параллелоэдров разбиения $\mathcal{T}(P),$ можно представить в виде суммы нескольких векторов из $\mathcal{N}(P).$

\begin{con}[Г.Вороной \cite{Vor}]\label{con2}
В множестве $\mathcal{N}(P)$ можно выбрать $d$ векторов, являющиеся базисом решетки $\Lambda.$
\end{con}

В данной работе мы докажем вторую гипотезу Вороного в частном случае зонотопов, являющихся параллелоэдрами.

\section{Разбиения Делоне и Вороного}

Построение многогранника Дирихле-Вороного из предыдущего раздела можно обобщить на случай произвольной квадратичной формы.

\begin{defin}
Пусть $\phi: \mathbb{R}^d \longrightarrow \mathbb{R}$ --- положительно определенная квадратичная форма в $\mathbb{R}^d$ и $O$ --- некоторая точка фиксированной решетки $\Lambda.$ {\it Многогранником Дирихле $P_\phi(\Lambda)$} для данной решетки $\Lambda$ относительно формы $\phi$ называется многогранник, состоящий из всех точек $X$ пространства $\mathbb{R}^d$, для которых значение формы $\phi$ на векторе $\overrightarrow{OX}$ не больше чем значение $\phi(\overrightarrow{O'X})$ для любой точки $O'\in\Lambda.$
\end{defin}

Построенный многогранник $P_\phi(\Lambda)$ также является параллелоэдром и соответствующее разбиение называется {\it разбиением Дирихле-Вороного} $\mathcal{V}_\phi(\Lambda)$. Обычный многогранник Дирихле-Вороного получится если в качестве формы $\phi$ рассмотреть обычную форму метрики $\phi(\ib x)=\ib x^T\ib x.$ При этом, если к данной решетке и данному многограннику $P_\phi(\Lambda)$ с формой $\phi(\ib x)=\ib x^T Q\ib x$ мы применим аффинное преобразование $\mathcal{A},$ то в результате получится решетка $\mathcal{A}\Lambda$ и многогранник $P_{\mathcal{A}^{-1}\phi}(\mathcal{A}\Lambda),$ где форма $\mathcal{A}^{-1}\phi$ задается равенством $\mathcal A^{-1}\phi(\ib x)=(\mathcal A^{-1}\ib x)^T Q (\mathcal A^{-1}\ib x).$ Следовательно, для того, чтобы доказать первую гипотезу Вороного для данного параллелоэдра $P$ с данной решеткой $\Lambda$ достаточно найти такую квадратичную форму $\phi,$ что $P=P_\phi(\Lambda).$

Двойственным к разбиению $\mathcal{V}_\phi(\Lambda)$ на многогранники Дирихле-Вороного $P_\phi(\Lambda)$ является разбиение Делоне, также определяемое по данным решетке $\Lambda$ и положительно определенной квадратичной форме $\phi$. Квадратичная форма $\phi$ определяет эллипсоид $\mathcal E_\phi$, задаваемый уравнением $\phi(\ib x)=1.$ Рассмотрим какой-нибудь эллипсоид гомотетичный $\mathcal E_\phi$ содержащий на границе $d$-мерное подмножество точек из $\Lambda,$ но не содержащий внутри себя ни одной точки $\Lambda.$ Такой эллипсоид определяет вписанный в него выпуклый многогранник. Множество таких многогранников, вписанных в ``пустые'' эллипсоид называется {\it разбиением Делоне} $\mathcal{D}_\phi(\Lambda).$

Если мы рассмотрим произвольную вершину $O$ разбиения Дирихле-Вороного, то значение квадратичной формы $\phi$ на векторах,  соединяющих $O$ с центрами многогранников из $\mathcal{V}_\phi(\Lambda)$ ее содержащих, равны, а значит $O$ и будет центром ``пустого'' эллипсоида, определяющего один из многогранников разбиения Делоне $\mathcal{D}_\phi(\Lambda).$

\section{Зонотопы и дайсинги}

\begin{defin}
Многогранник $P\subseteq \mathbb{R}^d$ называется {\it зонотопом} если его можно представить в виде суммы Минковского некоторого конечного набора отрезков. Зонотоп, являющийся суммой Минковского отрезков, образованных векторами множества $V=\{\ib v_1,\ldots,\ib v_n\},$ обозначается $Z(V).$ Векторы $\ib v_i$ называются {\it векторами зон}.

Эквивалентным образом зонотоп можно определить как проекцию куба $C^n$ некоторой размерности $n\geq d.$
\end{defin}

Эрдал в работе \cite{Erd} доказал первую гипотезу Вороного для зонотопов, являющихся параллелоэдрами. Позднее доказательство первой гипотезы Вороного для зонотопов было дано Деза и Гришухиным в терминах ориентированных матроидов \cite{DG}. В данной работе мы сформулируем и будем использовать некоторые факты о зонотопах и соответствующих разбиений Делоне и Дирихле-Вороного из указанной работы \cite{Erd}.

\begin{defin}
Рассмотрим в $d$-мерном пространстве $n$ семейств гиперплоскостей с нормалями $\ib d_1,\ldots, \ib d_n$ причем в каждом семействе все плоскости параллельны и дают разбиение пространства на слои равной толщины. Данный набор семейств и разбиение, которое они задают называются {\it дайсингом} $\mathfrak{D}(\ib d_1,\ldots,\ib d_n)$ если среди данных семейств можно выбрать $d$ с линейно независимыми нормалями и через любую точку, через которую проходят плоскости $d$ семейств, проходят плоскости и остальных семейств.
\end{defin}

Пусть точка с радиус-вектором $\ib a_0$ лежит на гиперплоскостях всех семейств дайсинга, тогда для данного дайсинга можно выбрать нормали $\ib d_1,\ldots,\ib d_n$ к гиперплоскостям его семейств таки образом, что все гиперплоскости дайсинга задаются уравнениями $\ib d_i\cdot(\ib x-\ib a_0)=a$ для всевозможных целых параметров $a.$ При этом каждый вектор $\ib d_i$ можно заменить на ему противоположный $- \ib d_i$.

\begin{defin}
В этом случае множество векторов $\mathcal{D}=\{\pm\ib d_1,\ldots,\pm\ib d_n\}$ называется {\it множеством нормалей} дайсинга $\mathfrak{D}(\ib d_1,\ldots,\ib d_n).$
\end{defin}

Выберем произвольную вершину дайсинга $\mathfrak{D}$ (то есть точку пересечения гиперплоскостей из $d$ семейств) и рассмотрим все ребра дайсинга, которые из нее выходят. Это множество ребер будет центрально симметричным.

\begin{defin}
Данное множество $\mathcal{E}=\{\pm\ib e_1,\ldots,\pm\ib e_k\}$ называется {\it множеством ребер} дайсинга. Очевидно, что построение этого множества не зависит от того какую мы вершину выбрали изначально.

Также отсюда следует, что все вершины дайсинга $\mathfrak{D}$ образуют решетку $\Lambda(\mathfrak{D}),$ которая называется {\it решеткой дайсинга}.
\end{defin}

В частности Эрдал доказал следующую теорему о связи множеств $\mathcal{D}$ и $\mathcal{E}$ \cite[Теорема 3.1]{Erd}.

\begin{thm}[Р. Эрдал]\label{edges}
Множество $\mathcal{D}$ с $d$ линейно независимыми векторами является множеством нормалей некоторого дайсинга если и только если существует такое множество $\mathcal{E}$, что:
\begin{enumerate}
\item[{\rm (E1)}] Любая пара противоположных векторов $\pm \ib e_i\in \mathcal{E}$ принадлежит одномерному пересечению $\ib d_{i_1}^\perp\cap\ldots\cap \ib d_{i_{d-1}}^\perp,$ где векторы $d_{i_j}\in \mathcal D$ и линейно независимы и наоборот для любых $d-1$ независимого вектора из $\mathcal{D}$ найдется соответствующая пара противоположных векторов из $\mathcal{E};$

\item[{\rm (E2)}] Для любых векторов $\ib d\in \mathcal D$ и $\ib e\in \mathcal E$ верно равенство $\ib d \cdot \ib e=0,\pm 1.$
\end{enumerate}
\end{thm}

\begin{defin}
Матрицы $\mathbf{D}$ и $\mathbf{E},$ столбцы которых состоят из векторов множеств $\mathcal D$ и $\mathcal{E}$ соответственно, называются матрицами {\it векторов нормалей} и {\it векторов ребер} дайсинга $\mathfrak{D}.$
\end{defin}

При аффинном преобразовании пространства $\mathbb{R}^d$ с матрицей $L$ матрицы $\mathbf{D}$ и $\mathbf{E}$ меняются по правилам $\mathbf{D}'=(L^{-1})^T\mathbf{D}$ и $\mathbf{E}'=L\mathbf{E}.$ Более того, существует аффинное преобразование после которого матрицы $\mathbf{D}$ и $\mathbf{E}$ состоят только из $0$ и $\pm1$ \cite[Теорема 3.3]{Erd}.

\begin{thm}[Р. Эрдал]
Существует аффинное преобразование, после которого матрица $\mathbf D'$ будет полностью унимодулярной, то есть любой ее минор будет равен $0$ или $\pm 1.$ Более того, это преобразование можно выбрать таким образом, что в каждом из множеств $\mathcal{D}$ и $\mathcal{E}$ можно выбрать по $d$ векторов, которые являются стандартным базисом пространства $\mathbb{R}^d,$ то есть векторами $(1,0,\ldots, 0)^T,(0,1,0,\ldots, 0)^T,\ldots,(0,\ldots,0,1)^T.$ Также в этом базисе элементами матриц $\mathbf{D}$ и $\mathbf{E}$ будут $0$ и $\pm 1.$
\end{thm}

Рассмотрим половину всех векторов из множества $\mathcal{D}$ среди которых нет противоположных. Данное множество $\mathcal{D}^+$ также однозначно определяет дайсинг $\mathfrak{D}$. Рассмотрим квадратичную форму
$$\phi(\ib x)=\sum_{\ib d\in \mathcal{D}^+}\omega_{\ib d}(\ib d\cdot \ib x)^2 $$ для некоторых положительных $\omega_{\ib d}.$ Тогда соответствующий многогранник Дирихле-Вороного описывается так \cite[Теорема 4.3]{Erd}:

\begin{thm}[Р. Эрдал]\label{zones}
Многогранником Дирихле-Вороного решетки $\Lambda(\mathfrak{D})$ относительно квадратичной формы $\phi(\ib x)$ является зонотоп с векторами зон $$\ib z_{\ib d}=\left(\sum_{\ib d\in\mathcal{D}^+}\omega_{\ib d}\ib d\ib d^T\right)^{-1}\omega_{\ib d}\ib d, \quad \ib d\in \mathcal D^+.$$
\end{thm}

Более того, в той же работе доказано, что если метрика $\phi(\ib x)$ --- стандартная евклидова, то векторы зон равны $\ib z_{\ib d}=\omega_{\ib d}\ib d$ (\cite[Разд. 6, стр. 442]{Erd}).

Также в \cite[Теорема 1.2]{Erd} доказана, что многогранник Дирихле-Вороного для решетки является зонотопом если и только если соответствующее разбиение Делоне является дайсингом.

\section{Вторая гипотеза Вороного для зонотопов}

\begin{lem}\label{face}
Пусть $Z=Z(V)$ --- $d$-мерный зонотоп. Любая гипергрань $Z$ порождается $(d-1)$-мерным подмножеством $U$ множества $V$ и наоборот, любой $(d-1)$-мерное подмножество $V$ порождает гипергрань $Z(U)$ зонотопа $Z(V).$
\end{lem}

\begin{proof}
Пусть $\pi$ --- гиперплоскость какой-либо гиперграни $Z.$
Зонотоп $Z(V)=Z(\ib v_1,\ldots \ib v_n)$ является проекцией куба $C^n\subset\mathbb{R}^n$ на пространство $\mathbb{R}^d$ вдоль $(n-d)$-мерного пространства $\psi.$ Рассмотрим в пространстве $\mathbb{R}^n$ гиперплоскость $\pi\times\psi;$ она является опорной гиперплоскостью для куба $C^n,$ а значит определяет его грань $F.$ Грань $F$ является кубом размерности $d-1$, который порождается сторонами куба $C^n,$ проецируемыми в векторы множества $V$, параллельные плоскости $\pi.$ Следовательно,  вся грань $F$ проецируется в параллельную копию зонотопа $Z(U)$ и в тоже время в грань зонотопа $Z(V).$ В обратную сторону доказательство аналогично.
\end{proof}

\begin{thm}
Для зонотопов, являющихся параллелоэдрами верна гипотеза $\ref{con2}.$
\end{thm}

\begin{proof}
Пусть зонотоп $Z$ является $d$-мерным параллелоэдром. Для него верна первая гипотеза Вороного \cite[Теорема 1.1]{Erd}, то есть существует такое аффинное преобразование $\mathcal{A},$ что зонотоп $\mathcal{A}Z$ является многогранником Дирихле-Вороного для некоторой решетки $\Lambda.$ Разбиение Делоне для решетки $\Lambda$ и стандартной метрики является дайсингом $\mathfrak{D}=\mathcal{D}(\pm\ib d_1,\ldots,\pm\ib d_n),$ следовательно (теорема \ref{zones}), векторы зон зонотопа $\mathcal{A}Z$ это векторы вида $\omega_i \ib d_i.$

Вектор, соединяющий центры $O_1$ и $O_2$ двух копий $\mathcal{A}Z$ с общей гипергранью $F$, перпендикулярен это гиперграни, так как все точки $F$ равноудалены от $O_1$ и $O_2.$ По лемме \ref{face} грани $F$ параллельны $d-1$ вектор из множества $\mathcal{D},$ а значит по теореме \ref{edges}, найдется вектор из $\mathcal{E}$ перпендикулярный $F,$  следовательно вектор $\overrightarrow{O_1O_2}$ лежит в $\mathcal E.$ Обратное также верно, если мы возьмем любой вектор из $\mathcal{E},$ то $d-1$ независимых векторов из $\mathcal{D},$ которые ему перпендикулярны, зададут гипергрань зонотопа $Z.$ Таким образом, множества $\mathcal{E}$ и $\mathcal{N}(Z)$ совпадают.

Рассмотрим унимодулярное представление соответствующего дайсинга. В данном представлении в множестве $\mathcal{E}$ содержится стандартный базис пространства, а остальные векторы множества $\mathcal{E}$ имеют целочисленные координаты. Следовательно, в унимодулярном представлении и в множестве $\mathcal{N}(Z)$ содержится стандартный базис пространства, а остальные векторы этого множества выражаются через данный базис с целыми коэффициентами.
\end{proof}